\documentclass[12pt]{article}
\title{Symplectic Resolutions for Nilpotent Orbits}
\author{Baohua FU}
\usepackage{amsmath,amssymb,amsthm,amscd}
\chardef\bslash=`\\
\newtheorem{Thm}{Theorem}[section]
\newtheorem{Cor}[Thm]{Corollary}
\newtheorem{Lem}[Thm]{Lemma}
\newtheorem{Prop}[Thm]{Proposition}
\newtheorem{Def}{Definition}[section]
\newtheorem{Rque}{Remark}

\def\cit{{\mathbb C}}
\def\nit{{\mathbb N}}
\def\qit{{\mathbb Q}}
\def\zit{{\mathbb Z}}
\def\pit{{\mathbb P}}
\def\0{{\mathcal O}}
\def\g{{\mathfrak g}}
\def\p{{\mathfrak p}}
\def\X{{\mathcal X}}
\def\C{{\mathcal C}}
\begin{document}
\maketitle
\begin{abstract}
In this paper, we prove that for any nilpotent orbit $\0$ in a semi-simple complex Lie algebra, equipped with the Kostant-Kirillov
symplectic form $\omega$, if for a proper resolution $\pi: Z \rightarrow \overline{\0}$, 
the 2-form $\pi^*(\omega)$ defined on $\pi^{-1}(\0)$ extends to a 
symplectic 2-form on $Z$, then $Z$ is isomorphic to the cotangent bundle $T^*(G/P)$ of a
 projective homogeneous space, and $\pi$ is the collapsing of the 
zero section.  Using this theorem, we determine all the 
varieties $\overline{\0}$  which admit such a resolution.  
We also calculate the Picard group of nilpotent orbits of classical type.  
\end{abstract}
\tableofcontents
\section*{Introduction}

Since A. Beauville's pioneering paper \cite{Be2}, symplectic singularities have received a particular attention by many mathematicians.
Recall that a holomorphic 2-form on a smooth variety is {\em symplectic} if it is closed and non-degenerate at every point. An algebraic variety 
$V$ (always over $k = \mathbb{C}$), smooth in codimension 1, is said to have {\em symplectic singularities} (or to be a {\em symplectic variety} or
a {\em symplectic singularity}) if there 
exists a holomorphic symplectic 2-form $\omega$ on $V_{reg}$ such that for any resolution of singularities $\pi: \widetilde{V} \rightarrow V$, the 2-form
$\pi^* \omega$ defined a priori on $\pi^{-1}(V_{reg})$ can be extended to a regular 2-form on $\widetilde{V}$. If furthermore 
the 2-form $\pi^* \omega$ extends to a symplectic 2-form on $\widetilde{V}$ for some resolution of $V$, then we say that $V$ admits a 
{\em symplectic resolution}.  

A resolution of singularities $\pi: \widetilde{V} \rightarrow V$ for a symplectic variety $V$ is called {\em crepant} if the canonical bundle of 
$\widetilde{V}$ is trivial. As we will show later, a resolution for a normal symplectic singularity
 is crepant if and only if it is symplectic. In particular, we see that
the existence of a symplectic resolution is independent of the special symplectic form on $V_{reg}$.

In general, it is difficult to determine whether a symplectic variety admits  a symplectic resolution or not. 
In \cite{CMSB}, Y. Cho, Y. Miyaoka and N. Shepherd-Barron have proved that any symplectic resolution for an isolated normal symplectic singularity
of dimension $2n \geq 4$ is analytically isomorphic to the resolution $T^*\pit^n \to \overline{\0}_{min}$, where $\overline{\0}_{min}$
is the closure of the minimal nilpotent orbit in $\mathfrak{sl}(n+1, \cit).$ It is important to find out more models for symplectic resolutions.

Examples of symplectic singularities are quotients of symplectic singularities by finite groups of automorphisms which preserve
 a symplectic 2-form on the regular locus. 
A particular case is the quotient of a complex vector space $\cit^{2n}$ by a finite group $G$ of symplectic automorphisms. The existence of 
symplectic resolutions for such quotient singularities has been studied by D. Kaledin in \cite{Ka1} and 
M. Verbitsky in \cite{Ver}. The motivation is to generalize the McKay correspondence
 to higher dimensions.  Strikingly, if $\cit^{2n}/G$ admits a symplectic
resolution, then the generalized McKay correspondence has been proved by D. Kaledin \cite{Ka2}.

Another important class of symplectic singularities is the closure $\overline{\0}$ of a nilpotent orbit $\0$ in a semi-simple
complex Lie algebra $\g$. As is well-known there exists a natural symplectic 2-form $\omega$ on any adjoint orbit $\0$ 
(identified with a co-adjoint orbit via the Killing form), 
called the Kostant-Kirillov form. 
 It was first proved by D. Panyushev (\cite{Pan}) that this symplectic form extends to any resolution, so $\overline{\0}$ is a 
variety with symplectic singularities. The main purpose of this paper is to find out all nilpotent orbit closures which admit a symplectic resolution.
 More precisely, we prove the following 
\begin{Thm}[Main theorem]
Let $\0$ be a nilpotent orbit in a semi-simple complex Lie algebra $\g$ and $G$ the adjoint group of $\g$. Suppose that we have a symplectic resolution
$\pi: Z \to \overline{\0}$. Then $Z$ is isomorphic to $T^*(G/P)$ for some parabolic subgroup $P$ of $G$, and  
under this isomorphism, the map $\pi$ corresponds to
$$T^*(G/P) \simeq G \times^P \mathfrak{u} \rightarrow \g, \qquad (g, X) \mapsto Ad(g)X ,$$ where $\mathfrak{u}$ is the nilradical of $\mathfrak{p}=Lie(P)$.
\end{Thm}

Recall that an element $X \in \mathfrak{u}$ is a {\em Richardson element} if $dim(G\cdot X) = 2 dim(G/P)$, or equivalently if $P\cdot X$ is dense
 in $\mathfrak{u}$ and $P$  contains the identity component of the centralizer $Z_G(X)$ (see \cite{Hes}). 
 The orbit $G\cdot X$ is called a {\em Richardson orbit}, which plays an important role in representation theory.
 It was R. Richardson who has shown that every
parabolic subgroup in $G$ has Richardson elements (see also  Proposition \ref{symI-km}).
As a direct corollary,  our theorem implies that if the closure of a nilpotent orbit $\0$ admits a symplectic resolution, then $\0$ is a Richardson orbit. 

The key point of the proof is to study the $\cit^*$-action on $\overline{\0}$. Using the same idea as in \cite{Ka1} and \cite{Ver},
we show that this action can be lifted to any symplectic resolution $Z$ of $\overline{\0}$. Then we use some standard analysis for this
action on $Z$, as done in \cite{Ka1} and \cite{Na2}.

Let us give a brief outline of the contents of this paper.

$\bullet$ Section 1 recalls some basic definitions. We  prove some easy propositions which reveal the relationship between crepant resolutions, symplectic 
resolutions and resolutions with small exceptional set. 

$\bullet$ In section 2 we calculate Picard groups of nilpotent orbits of classical type. 
 Explicit formulas are given. Then we use these results to determine when the
 normalization $\widetilde{\0}$ is $\qit$-factorial or factorial. 

$\bullet$ In section 3, firstly we prove the main theorem. Then we use it and some results of W. Hesselink \cite{Hes} to determine, in terms
of the partition corresponding to the orbit $\0$, all nilpotent orbit closures   $\overline{\0}$ (not of the type $E_7$ and $E_8$) 
which admit a symplectic resolution. \\

{\em Acknowledgments:} I want to thank A. Hirschowitz, J. Kock, C. Margerin and C. Pauly for helpful discussions. I am especially grateful to 
A. Beauville for many  valuable discussions and suggestions. Without his help, this work could never have been done.
I want to thank the referee for some pertinent remarks.  

In the printed version (Invent. Math. 151 (2003), 167-186) of this paper, there was a gap in the proof of Lemma \ref{gap}, which is pointed
out to me  by  Y. Namikawa. We add here Proposition \ref{orbitcover} (the idea of the proof is due to M. Brion) 
to fix this gap. With the suggestions of
M. Brion and Y. Namikawa, the final step in our previous proof is shortened with the help of a result in [BK]. 
Some errors in the calculs for the Picard group of a classical nilpotent orbit are also corrected. 

\section{Preliminaries}

Let $V$ be an irreducible complex algebraic variety. A morphism $\pi: \widetilde{V} \rightarrow V$ is called a {\em resolution} 
(or {\em desingularization})
 if $\pi$ is proper and $\pi$ induces an isomorphism outside the singular locus $Sing(V)$ of $V$. Since we are working over $\cit$, Hironaka's big theorem 
says that every $V$ admits a  desingularization. In this paper we will consider a particular class of resolutions 
(called symplectic resolutions). Suppose that 
the canonical divisor $K_V$ of $V$ is a Cartier divisor. In this case, there is the following notion.
\begin{Def}
A resolution $\pi: \widetilde{V} \rightarrow V$ is called {\em crepant} if $\pi^*(K_V) = K_{\widetilde{V}}$, i.e. if $\pi$ preserves the canonical class.
\end{Def}

 Recall that a holomorphic 2-form on a smooth variety is {\em symplectic} if it is closed and non-degenerate at every point. 
\begin{Def}
An algebraic variety $V$, smooth in codimesion 1,
 is said to have {\em symplectic singularities} (or to be a {\em symplectic variety} or a {\em symplectic singularity}) if there exists a holomorphic
symplectic 2-form  $\omega$ on $V_{reg}$ such that for any resolution 
$\pi: \widetilde{V} \rightarrow V$, the 2-form $\pi^*\omega$ defined a priori on $\pi^{-1}(V_{reg})$ extends to 
a holomorphic 2-form $\Omega$ on $\widetilde{V}$.
\end{Def}
\begin{Rque}
In \cite{Be2}, one requires that a symplectic singularity is normal. However, since we deal essentially with 
$\overline{\0}$ in this paper, which is in general not normal (see \cite{KP1}, \cite{KP2}), we decide to adopt this weaker notion
of symplectic singularities.
\end{Rque}

Basic examples of varieties with symplectic singularities are quotients of symplectic singularities by finite groups of automorphisms, preserving a 
symplectic 2-form on the regular locus(see \cite{Be2}). Some classification theorems for isolated symplectic singularities have been
proved by S. Druel since the original paper \cite{Be2}.
\begin{Def}
A resolution  $\pi: \widetilde{V} \rightarrow V$ for a symplectic variety $(V, \omega)$ is called {\em symplectic}
 if the 2-form $\Omega$ extending $\pi^*\omega$
is a symplectic 2-form on $\widetilde{V}$.
\end{Def}

Note that for a normal symplectic variety $V$, its canonical sheaf is trivial, thus a symplectic resolution is a crepant resolution.
The converse is also true, as shown by the following proposition (see also \cite{Ka1} and \cite{Ver}).
\begin{Prop}\label{sym=crep}
Let $(V, \omega)$ be a normal symplectic variety of dimension $2n$ and $\pi: \widetilde{V} \rightarrow V$ a resolution,  
 then $\pi$ is a crepant resolution if and only if $\pi$ is a symplectic resolution.
\end{Prop}
\proof
We have seen that a symplectic resolution is a crepant resolution. 
Now suppose that $\pi$ is crepant.  Since $\omega^n$ has no zeros on $V_{reg}$, it extends to a global section of $K_V$.
That $\pi$ is crepant implies that $\pi^*(\omega^n)$ extends to a global section on $\widetilde{V}$ without zeroes. By our assumption, $V$ has symplectic 
singularities, so  $\pi^*\omega$ extends to a 2-form $\Omega$ on $\widetilde{V}$, furthermore we have $\Omega^n =\pi^*(\omega^n)$, thus $ \Omega^n$ has no 
zeroes on $\widetilde{V}$, which implies that $\Omega$ is non-degenerate everywhere. 
Now on $\pi^{-1}(V_{reg})$, we have $d \Omega = d \pi^*(\omega) = \pi^*(d\omega) = 0$. Since $\pi^{-1}(V_{reg})$ is open-dense,  $d \Omega =0$ on
 $\widetilde{V}$, i.e. the 2-form $\Omega$  is closed, which gives that $\Omega$ is a symplectic 2-form, i.e. $\pi$ is a symplectic resolution. 
\qed \\

In particular, we see that the existence of a symplectic resolution is independent of the special symplectic form on $V_{reg}$. 
This easy proposition gives the following criterion for a resolution being symplectic in the case of $Sing(V)$ having higher codimension.
\begin{Prop}
Let $V$ be a normal  symplectic variety with   $codim(Sing(V)) \geq 3$. Then $V$ has a symplectic resolution if and only if 
there exists a resolution  $\pi: \widetilde{V} \rightarrow V$ such that $codim(\pi^{-1}(Sing(V))) \geq 2$.
\end{Prop}
\proof
Let $\pi: \widetilde{V} \rightarrow V$ be a symplecitc resolution. Suppose that the exceptional set $E = \pi^{-1}(Sing(V))$ has a component of codimension 1,
 then by corollary 1.5 of \cite{Nam} (where the proof works without the assumption ``projective'') the image of this component by $\pi$ in $V$ should be of
codimension 2. This is impossible since $\pi$ is a resolution of singularities and $codim(Sing(V)) \geq 3$, so we have $codim(E) \geq 2$.

Conversely, if $\pi$ is resolution such that $codim(\pi^{-1}(Sing(V))) \geq 2$, then $\widetilde{V} - \pi^{-1}(V_{reg})$ has codimension $\geq 2$, so 
$K_{\widetilde{V}} = K_{\pi^{-1}(V_{reg})} = \pi^*(K_{V_{reg}})$, which shows that $K_{\widetilde{V}}$ is trivial. 
Thus the resolution $\pi$ is crepant, which
is also symplectic by the above proposition.
\qed \\

Recall that a variety is {\em $\qit$-factorial} if any Weil divisor has a multiple that is a Cartier divisor. A smooth variety is $\qit$-factorial.
\begin{Cor}\label{symI-cor}
A $\qit$-factorial normal symplectic variety $V$ whose singular locus $Sing(V)$ is not of pure co-dimension 2
does not admit any symplectic resolution.  
\end{Cor}
\proof
By an argument of \cite{Deb} (p.28), $V$ being normal and $\qit$-factorial implies that for any birational map $\pi: Z \rightarrow V$, the 
exceptional set $E \subset Z$ of $\pi$ is of pure codimension 1 in $Z$. Now the corollary follows in the same way as the preceding proposition.  
\qed \\

Here we recall some basic notions concerning nilpotent orbits. A detailed and excellent discussion can be found in \cite{CM}. Let $\g$ be a 
semi-simple complex 
Lie algebra and let $G_{ad}$ (or $G$ for short) be the identity component of its automorphism group, which is called the adjoint group of $\g$.
Recall that each adjoint orbit $\0$ (identified with a co-adjoint orbit in $\g^*$ via the Killing form) carries a canonical $G$-invariant 
symplectic 2-form $\omega$, the Kostant-Kirillov  form. An adjoint orbit is closed
if and only if it is a semi-simple orbit. For a nilpotent orbit $\0$, its closure $\overline{\0}$ in $\g$ is not necessarily normal (see \cite{KP1},
\cite{KP2}). 
It is shown in \cite{Pan} (see also \cite{Be1}) that the closure  $\overline{\0}$ of $\0$ is a symplectic variety.

Another feature of nilpotent orbits is the existence of a $\cit^*$-action. This can be seen as follows. Take a nilpotent element $X \in \g$, 
we want to show that $\lambda X$ is conjugate to $X$ for any $\lambda \in \cit^*$.
By the Jacobson-Morozov theorem, there exists a standard triple $(H,X,Y)$ for $X$, i.e. we have $$[H,X] = 2X, \quad [H,Y] = -2Y, 
\quad [X,Y] = H. $$  So we have an isomorphism $\phi: \mathfrak{sl}_2 \rightarrow \cit\langle H,X,Y \rangle$. 
Now it is clear that $\lambda X = gXg^{-1}$ where  $g = exp(cH)$ with $c \in \cit $ satisfying $exp(2c) = \lambda$.
For this action, we have 
\begin{Lem}
For $\lambda \in \cit^*$, we have $\lambda^* \omega = \lambda \omega$.
\end{Lem}
\proof
For an element $X \in \g$, let $\xi^X$ be the vector field $$\xi^X(Z) = \frac{d}{dt}|_{t=0} exp(tX) \cdot Z = [X,Z].$$ Then 
for any $\lambda \in \cit^*$, we have $\lambda_*\xi^X(Z) = \frac{d}{dt}|_{t=0} \lambda exp(tX) \cdot Z = [X, \lambda Z] = \xi^X(\lambda Z)$.
Let $\kappa(\cdot, \cdot)$ be the Killing form on $\g$, then by definition $$\omega_Z(\xi^X(Z), \xi^Y(Z)) = \kappa(Z,[X,Y]).$$
Now $$
\begin{aligned}
\lambda^* & \omega_Z(\xi^X(Z),\xi^Y(Z))  = \omega_{\lambda Z}(\lambda_*(\xi^X(Z)), \lambda_*(\xi^Y(Z))) \\ & = \omega_{\lambda Z}(\xi^X(\lambda Z),
 \xi^Y(\lambda Z))  = \kappa(\lambda Z, [X,Y]) = \lambda\omega_Z(\xi^X(Z), \xi^Y(Z)), \end{aligned}$$
so we have  $\lambda^* \omega = \lambda \omega$.
\qed \\

Now $\g$ is an $\mathfrak{sl}_2$-module via  the above isomorphism $\phi$, so  $\g$ is decomposed as $\g = \oplus_{i \in \zit} \g_i$, 
where $\g_i = \{ Z \in \g |\, [H,Z] = i Z \}.$  The nilpotent orbit $\0$ is called {\em even} if $\g_1=0$, or equivalently if $\g_{2k+1} = 0$ for
all $k \in \zit$ (see lemma 3.8.7 \cite{CM}). 
Let $\mathfrak{p} = \oplus_{i \geq 0} \g_i$ and $P$ a connected subgroup in $G$ with Lie algebra $\mathfrak{p}$. Then $\mathfrak{p}$
is a parabolic sub-algebra of $\g$ and $P$ is a parabolic subgroup of $G$.
Set $\mathfrak{n} = \oplus_{i \geq 2} \g_i$.
Then for an even orbit $\0$, there is an isomorphism between $(\g/\mathfrak{p})^*$ and $\mathfrak{n}$, and the latter
is nothing but the nilradical of $\mathfrak{p}$.
\begin{Prop}[Springer's resolution]
Let $\0$ be an even nilpotent orbit in $\g$, then there exists a $G$-equivariant resolution of singularities 
$$\pi: T^*(G/P) \simeq G\times^P\mathfrak{n} \rightarrow \overline{\0}, \quad (g,X) \mapsto Ad(g)X.$$
\end{Prop}

Thus the variety $\overline{\0}$ admits a symplectic resolution. 
The motivation of this paper is to find out all 
varieties $\overline{\0}$ which admit a symplectic resolution. This is achieved at the end of section 3. 

\section{Picard groups and $\qit$-factority}
\subsection{Picard groups of nilpotent orbits}
 Let $\g$ be a complex simple Lie algebra, and $G$ the adjoint group of $\g$.
For a nilpotent element  $X \in \g$, we denote by $\0_X$ the nilpotent orbit $G \cdot X $. It is isomorphic to $G/G^X$, 
where $G^X$ is the stabilizer of $X$ in $G$. 
The purpose of this section is to calculate Picard groups for these nilpotent orbits in case of $\g$ being of classical type.

Instead of working with $G$, we will work with the universal covering $G_{sc}$ of $G$. In this case we have 
$\0_X = G_{sc}/G_{sc}^X$, where $G_{sc}^X$ is the stabilizer of $X$ in $G_{sc}$. 

Recall that  $\g$ is decomposed as $\g = \oplus_{i \in \zit} \g_i$, where $\g_i = \{ Z \in \g |\, [H,Z] = i Z \}.$
Set $\mathfrak{u}^X = \oplus_{i > 0} \g_i^X$, where $\g_i^X = \g^X \cap \g_i$, and set 
$\g^{\phi} = \{ Z \in \g |\, [Z,X] = [Z,Y] = [Z,H] = 0 \}$, then $\g^X = \mathfrak{u}^X \oplus \g^\phi$.  Let $U^X$ be the connected subgroup of 
$G_{sc}^X$ with Lie algebra $\mathfrak{u}^X$, and let $G_{sc}^\phi$ be the centralizer of $Im(\phi)$ in $G_{sc}$, then we have 
\begin{Prop}[Barbasch-Vogan, Kostant]
There is a semi-direct product decomposition $G_{sc}^X = U^X \cdot G_{sc}^\phi$.
\end{Prop}

Recall that the character group $\X(G)$ of an algebraic group $G$ is  the abelian group of algebraic group morphisms
 between $G$ and $\cit^*$, i.e. $\X(G) = Hom(G, \cit^*)$. 
\begin{Lem}\label{symreso-lem22}
We have $\X(G_{sc}^X) = \X(G_{sc}^\phi)$.
\end{Lem}
\proof
 Note that the algebra $\mathfrak{u}^X$ is nilpotent, so the group $U^X$ is unipotent. As is well-known, every unipotent group has trivial character group.
Now our lemma follows from the above proposition.
\qed \\

So to calculate $\X(G_{sc}^X)$, we need to calculate the character group $\X(G_{sc}^\phi)$. In the classical cases, we can describe explicitly 
the subgroup $G_{sc}^\phi$. Before giving the description, we need some notations. For any group $H$, 
let $H_\Delta^m$ denote the diagonal copy of $H$ in $H^m$. If 
$H_1, \cdots, H_m$ are matrix groups, let $S(\prod_i H_i)$ be the subgroup of $\prod_i H_i$ consisting of $m$-tuples of matrices whose determinants have 
product 1.

Recall that a {\em partition} {\bf d} of $n$ is a tuple of integers
$[d_1, \cdots, d_N]$ such that $d_1 \geq d_2 \geq \cdots \geq d_N > 0$ and $\sum_{j=1}^N d_j = n$.
In classical cases, a nilpotent orbit can be encoded by some partition {\bf d}. To illustrate the idea, let us consider the case of $\mathfrak{sl}_n$. 
Every nilpotent element $X \in \mathfrak{sl}_n$ is conjugate to an element of the form $diag(J_{d_1}, \cdots, J_{d_N})$, where $J_{d_j}$ is a Jordan block
of type $d_j \times d_j$ and {\bf d}= $[d_1, \cdots, d_N]$ is a partition of $n$.  It is clear that this partition is invariant under conjugation,
thus to the nilpotent orbit $\0_X$ we associate the partition {\bf d}, which establishes a bijection between nilpotent orbits in $\mathfrak{sl}_n$
and partitions of $n$. A similar bijection exists for nilpotent orbits in $\mathfrak{sp}_{2n}$ and $\mathfrak{so}_m$ (see section 5.1 of \cite{CM}).
It is easy to see that the orbit $\0_X$ is even if and only if all the parts $d_i$ have the same parity.
For a partition {\bf d} = $[d_1, \cdots, d_N]$, we put $r_i = \# \{j | d_j = i\}$ and $s_i = \# \{j | d_j \geq i \}$. Then we have 
(see \cite{CM} theorem 6.1.3): 
\begin{Prop}[Springer-Steinberg]
\begin{equation*} G_{sc}^\phi = \begin{cases} S(\prod_i (GL_{r_i})_\Delta^i)  & \g  = \mathfrak{sl}_n; \\
 \prod_{i \, odd} (Sp_{r_i})^i_\Delta \times \prod_{i\,even}(O_{r_i})^i_\Delta   & \g =\mathfrak{sp}_{2n}; \\
 double\ cover \ of\  S(\prod_{i \, even} (Sp_{r_i})^i_\Delta \times \prod_{i\,odd}(O_{r_i})^i_\Delta)  & \g = \mathfrak{so}_{m}. 
\end{cases} \end{equation*}
\end{Prop}
\begin{Lem}
If $\g = \mathfrak{sl}_n$, then $\X(G_{sc}^X) = \zit^{k-1} \oplus \zit/ c \zit $, where $c= gcd(d_1, \cdots, d_N)$ and 
$k = \# \{i | r_i \neq 0 \}$ is the number of distinct $d_i$ in {\bf d}.
\end{Lem}
\proof
Recall that $\X(GL_r) \cong \zit$ for any $r > 0$, where the isomorphism is given by 
$\zit \ni l \mapsto (A \mapsto det(A)^l)$. So $\X(\prod_i (GL_{r_i})_\Delta^i) \cong \zit^{k}$, where $k = \# \{i | r_i \neq 0 \}$ is the number of
 distinct $d_i$ in the partition {\bf d}. Now let $G^\circ$ be the identity component of $S(\prod_i (GL_{r_i})_\Delta^i)$, then 
$\X(G^\circ) = \zit^{k-1}$. By 6.1.2(\cite{CM}), $G_{sc}^\phi/G^\circ$ is nothing but $\pi_1 (\0_{\bf d}) = \zit/ c \zit.$
Notice that $\X( \zit/ c \zit) = \zit/ c \zit$, which gives the result. 
\qed \\
\begin{Prop}\label{symI-pro1}
Let $\0_X$ be a nilpotent orbit in $\mathfrak{sl}_n$ corresponding to the partition {\bf d}. Let $c= gcd(d_1, \cdots, d_N)$ and
 $k$ the number of distinct $d_i$ in {\bf d}.
Then we have $Pic(\0_X) = \zit^{k-1} \oplus \zit/ c \zit $.
\end{Prop}
\proof
By a result of V. Popov (\cite{Pop}), the Picard group of a connected semisimple Lie group is isomorphic to its fundamental group. So the Picard group
$Pic(G_{sc})$ of the simply connected semi-simple group $G_{sc}$ is trivial.
Now  the following exact sequence (\cite{KKV}): 
$$ 0 \rightarrow \X(G_{sc}^X) \rightarrow Pic(\0_X) \rightarrow  Pic(G_{sc}) = 0 ,$$
gives that $Pic(\0_X) = \X(G_{sc}^X)$, which is equal to $\zit^{k-1}\oplus \zit/ c \zit$ by the above lemma.
\qed \\

Now we consider the case of a simple Lie algebra of B-C-D type. 
\begin{Thm}
Let $\g$ be a simple complex Lie algebra  of B-C-D type, and $\0_X$ a nilpotent orbit in $\g$, then we have  
$$Pic(\0_X) = \X(\pi_1(\0_X)) \oplus \zit^l,$$ where $\pi_1(\0_X)$ is the fundamental group of $\0_X$ and 
$l$ is the number of even (resp. odd) parts which appears exactly twice in {\bf d} for   $\g =\mathfrak{sp}_{2n}$
(resp. $\g = \mathfrak{so}_m$). 
\end{Thm}
\proof   
By our lemma \ref{symreso-lem22}, we have $\X((G_{sc}^X)^\circ) = \X((G_{sc}^\phi)^\circ) $. By the above theorem of Springer-Steinberg, we see that
$(G_{sc}^\phi)^\circ$ is product of copies of semi-simple Lie groups with  $l$ copies of $SO_2 \simeq \cit^*$,
 thus $\X((G_{sc}^\phi)^\circ) = \zit^l$. 

By lemma 6.1.1 of \cite{CM}, the component group $G_{sc}^X / (G_{sc}^X)^\circ$ is isomorphic to the fundamental group
$\pi_1(\0_X)$ of $\0_X$, thus we have $\X(G_{sc}^X / (G_{sc}^X)^\circ) = \X(\pi_1(\0_X))$. Now a similar argument as we did in the above proof
gives that $Pic(\0_X) = \X(G_{sc}^X) = \X(\pi_1(\0_X)) \oplus \zit^l$.
\qed \\

Combining with explicit formulas for $\pi_1(\0_X)$ given in \cite{CM} (Cor. 6.1.6), we can give the following  formulas for
$Pic(\0_X)$. Recall that we say a partition {\bf d} is {\em rather odd} if all of its odd parts have multiplicity one. 
\begin{Cor}\label{symI-cor2}
Let $\0_X$ be a nilpotent orbit in $\g$ corresponding to the partition {\bf d}. Let $a$ be the number of distinct odd $d_i$ 
and let $b$ be the number of distinct even $d_i$. Let $l$ be the number of even (resp. odd) parts which appears exactly twice in 
{\bf d} for   $\g =\mathfrak{sp}_{2n}$ (resp. $\g = \mathfrak{so}_m$). 
Then we have  

(1) For $\g = \mathfrak{sp}_{2n}$, we have $Pic(\0_X) = (\zit/2\zit)^b \oplus \zit^l$;

(2) For $\g = \mathfrak{so}_{2n+1}$, if {\bf d} is rather odd, then $Pic(\0_X)$ is an extension of $\zit/2\zit$ by $(\zit/2\zit)^{a-1}$,
 otherwise $Pic(\0_X) = (\zit/2\zit)^{a-1}\oplus \zit^l$;  

(3) For $\g = \mathfrak{so}_{2n}$, if {\bf d} is rather odd, then $Pic(\0_X)$ is an extension of $\zit/2\zit$ by $(\zit/2\zit)^{max\{0,a-1\}}$,
otherwise $Pic(\0_X) = (\zit/2\zit)^{max\{0,a-1\}}\oplus \zit^l$.
\end{Cor}

\subsection{$\qit$-factority and symplectic resolutions}

Let $\g$ be a simple Lie algebra of classical type, and $\0_X$ be a nilpotent orbit. The closure $\overline{\0}_X$ is not always normal 
(for precise results see \cite{KP}). We denote by
$\widetilde{\0}_X$ its normalization. 
In this section, we will give some results on the properties of being $\qit$-factorial and factorial for the normal variety $\widetilde{\0}_X$. 

For an irreducible algebraic variety $V$, which is smooth in codimension 1, we denote by $Cl(V)$ its divisor class group, i.e. the Weil divisors modulo linear
equivalences (\cite{Har}). The following lemma is easily proved.
\begin{Lem}
Let $V$ be an irreducible algebraic variety, which is smooth in codimension 1. Let $\pi: \widetilde{V} \rightarrow V$ be the normalization. 
Then the induced map $\pi^*: Cl(V) \rightarrow Cl(\widetilde{V})$ is an isomorphism of groups.
\end{Lem}
\begin{Prop}\label{symI-pro2}
Let $\g$ be a simple Lie algebra of B-C-D type, and let $\0_X$ be a nilpotent orbit in $\g$
such that $l=0$, then the normal variety $\widetilde{\0}_X$ is $\qit$-factorial.
\end{Prop}
\proof
Since la codimension of $\overline{\0} - \0$ is at least 2, $\overline{\0}_X$ is smooth in codimension 1 (thus $Cl(\overline{\0}) = Pic(\0)$).
By  the preceding lemma, we have $Cl(\widetilde{\0}_X) \cong Cl(\overline{\0}_X)$. Now by  
proposition II.6.5 of \cite{Har}, we have $Cl(\overline{\0}_X) \cong Cl(\0_X)$. That $\0_X$ is smooth gives  $Cl(\0_X) = Pic(\0_X)$, thus 
$Cl(\widetilde{\0}_X) \cong Pic(\0_X)$. Now by our calculations in the preceding section, we know that $Pic(\0_X)$ is a finite group, so 
$Cl(\widetilde{\0}_X) \otimes_\zit \qit  \cong Pic(\0_X)\otimes_\zit \qit = 0$. In particular, the map 
$$Pic(\widetilde{\0}_X)\otimes_\zit \qit \rightarrow Cl(\widetilde{\0}_X) \otimes_\zit \qit$$
is surjective, hence  $\widetilde{\0}_X$ is $\qit$-factorial. 
\qed \\
\begin{Rque}
The condition that $l=0$ cannot be dropped, as shown by the following example (communicated to me by Y. Namikawa):
Let $\0$ be the nilpotent orbit in $\mathfrak{so}_{2n}$ corresponding to the partition $[2^{n-1}, 1^2]$,
where $n\geq 3$ is an odd integer. Then 
$\widetilde{\0}$ admits a small  symplectic resolution $T^* G_{iso}^+(n, 2n) \to \widetilde{\0}$, where
$G_{iso}^+(n, 2n)$ is a connected component of the orthogonal Grassmannian $G_{iso}(n, 2n)$.  Thus $\widetilde{\0}$
is not $\qit$-factorial. 
\end{Rque}

\begin{Prop}
Let $\g$ be a simple Lie algebra, and $\0_X$ be a non-zero nilpotent orbit in $\g$ corresponding to the partition {\bf d} = $[d_1, \cdots, d_N]$. Then we have:

(1) For $\g = \mathfrak{sl}_n$, the normal variety $\widetilde{\0}_X$ is never factorial. It is $\qit$-factorial 
if $d_1=d_2=\cdots=d_N$, i.e. there is only one distinct $d_i$;

(2) For $\g = \mathfrak{sp}_{2n}$, the normal variety  $\widetilde{\0}_X$ is factorial if and only if every $d_i$ is odd;

(3) For $\g = \mathfrak{so}_{2n} $, the normal variety  $\widetilde{\0}_X$ is factorial if and only if there exists exact one distinct odd $d_i$
with multiplicity at least 4;

(4) For $\g = \mathfrak{so}_{2n+1}$, the normal variety  $\widetilde{\0}_X$ is factorial if and only if there exists just one distinct odd $d_i$ with 
multiplicity at least 3.
\end{Prop}
\proof
From proposition II.6.2 of \cite{Har}, the affine normal variety $\widetilde{\0}_X$ is factorial if and only if $Cl(\widetilde{\0}_X) = 0 $.
By the proof of the above proposition, this is equivalent to $Pic(\0_X) = 0$. Now we just do a case-by-case check, based on our proposition \ref{symI-pro1} and 
corollary \ref{symI-cor2}. For example, when $\g = \mathfrak{so}_{2n+1}$, $Pic(\0_X)=0$ if and only if {\bf d} is not rather odd and $a = 1, l=0$,
 i.e. {\bf d} has only
one distinct odd $d_i$, with multiplicity at least 2. But the sum $\sum d_i = 2n+1$ is odd, so the multiplicity should be at least 3. 
\qed \\

 Applying these results to symplectic resolutions, we have
\begin{Prop}
Let $\g$ be a simple Lie algebra of type B-C-D, and $\0_X$ be a nilpotent orbit in $\g$ with $l=0$.
 If $\overline{\0}_X - \0_X$ is not of pure co-dimension 2
in  $\overline{\0}_X$, then the  normal variety $\widetilde{\0}_X$  does not admit any symplectic resolution. 
\end{Prop}
\proof
This comes directly from the above proposition \ref{symI-pro2} and corollary \ref{symI-cor}.
\qed \\
{\em Examples:}  In the case of $\g = \mathfrak{so}_8$, we find 2 orbits which do not admit any symplectic resolution, corresponding to 
the partitions $[3,2^2,1]$ and $[2^2,1^4]$, while the other nilpotent orbits are all even, 
thus having a symplectic resolution by Springer's resolution.  

Recall the the minimal nilpotent orbit $\0_{min}$ in $\mathfrak{so}_{2n+1} (n \geq 2)$ (resp. $\mathfrak{sp}_{2n} (n \geq 3)$, 
$\mathfrak{so}_{2n} (n \geq 4)$)
is the nilpotent orbit corresponding to the partition $[2^2, 1^{2n-3}]$ (resp. $[2, 1^{2n-1}]$, $[2^2, 1^{2n-4}]$). Thus the precedent
proposition implies:
\begin{Cor}
Let $\g$ be a simple Lie algebra of type B-C-D. The closure of the minimal orbit $\overline{\0}_{min}$ in $\g$ does not admit any symplectic resolution.
\end{Cor}
\begin{Rque}
This is also true for minimal orbits in exceptional Lie algebras by a similar argument. We just need to note that the Picard groups for these cases are
 trivial. This follows also from our discussions in section 3.4.   
\end{Rque}

As we will see later, the closure of every nilpotent orbit in $\mathfrak{sl}_n$ admits a symplectic resolution. This is more or less 
known to some experts.

\section{Symplectic resolutions for nilpotent orbits}
\subsection{Lifting the action of $\cit^*$ and $G$}
Let $\g$ be a semi-simple complex Lie algebra and $G$ its adjoint group. Let $\0$ be a nilpotent orbit in $\g$ and $\overline{\0}$ its Zariski closure in $\g$.
   Recall that there exists an action of $G$ (resp. $\cit^*$) on $\overline{\0}$. The purpose of this section is to prove 
that this action can be 
lifted to any symplectic resolution of $\overline{\0}$. Let $\pi: Z \rightarrow \overline{\0}$ be a symplectic resolution, and
$\Omega$ the symplectic 2-form on $Z$ extending $\pi^*\omega$, where $\omega$ is the Kostant-Kirillov symplectic form on $\0$.
\begin{Prop} \label{symI-lift}
The action of $G$ (resp. $\cit^*$) on $\overline{\0}$ lifts to $Z$, in such a way that  $\pi$ is $G$-equivariant (resp. $\cit^*$-equivariant). 
\end{Prop}
\proof
The proof is inspired from the proof of theorem 1.3 in \cite{Ka1} and theorem 2.5 in \cite{Ver}. Firstly we will show that the action of $\g$ 
on $\overline{\0}$ can be lifted to an action
of $\g$ on $Z$. To this end, let $X \in \g$. Consider the vector field $\xi^X$ on $\0$ defined by $\xi^X(Y) = \frac{d}{dt}|_{t=0} exp(tX)\cdot Y = [X,Y]$,
 where $\cdot$
denotes the adjoint action of $G$ on $\g$. Now the symplectic form $\omega$ on $\0$ gives an isomorphism $\Omega^1(\0) \simeq \mathcal{T}(\0)$ between
1-forms and vector fields on $\0$. 
Let us denote by
$\alpha_X$ the 1-form on $\0$ corresponding to the vector field $\xi^X$. The key point is the following claim. \\

{\em Claim:} the 1-form $\pi^*(\alpha_X)$ extends to the whole of $Z$.

\proof[Proof of the claim] Take a Hermitian metric $h$ on $Z$. Set $U = \pi^{-1}(\0)$. Then $\pi^*(\alpha_X)$ can be extended to the whole of 
$Z$ unless it has singularities on the complement $Z - U$. So we need to show that for any compact set $K \subset Z$ and for any $z \in U \cap K$, 
the Hermitian norm of $\pi^*(\alpha_X)$ is bounded by some constant depending on $K$. 

Since $\pi$ is analytic, thus Lipschitz on compact subsets $K \subset Z$. By rescaling the metric $h$, we can suppose that $\pi|_K$ decreases distance. 
Here the metric on $\0$ is the one induced from the metric on $\g$. We need to show that on $\pi(K)\cap \0$, the 1-form $\alpha_X$
has bounded norm.  This follows from the fact that the vector field $\xi^X$ has bounded norm on $\pi(K) \cap \0$ and $\omega$ is $G$-invariant, so it 
 also has bounded norm on $\pi(K) \cap U$. 
\qed \\

Let us denote by $\tilde{\alpha}$ the extended 1-form on $Z$. The symplectic form $\Omega$ on $Z$ gives a 1-1 correspondence between 1-forms and 
vector fields on $Z$,
so we get a vector field $\zeta^X$ on $Z$. 
This gives an action of $G_{sc}$ on $Z$. The center of $G_{sc}$ acts trivially on the open-dense set $\pi^{-1}(\0)$, so it acts trivially on the whole of $Z$.
 Since $G$ is the quotient of  $G_{sc}$ by its center, we get an action of $G$ on $Z$. 

A similar argument shows that we can also lift  the $\cit^*$-action to $Z$.  It is clear from the construction that $\pi$ is
 $G$-equivariant (resp. $\cit^*$-equivariant). 
\qed \\
\begin{Cor}
The symplectic form $\Omega$ is $G$-invariant, and for the $\cit^*$-action, we have $\lambda^* \Omega = \lambda \Omega$ for any $\lambda \in \cit^*$.
\end{Cor}
\proof
This comes from the above proposition and the corresponding properties of $\omega$, which is $G$-invariant and 
satisfies $\lambda^* \omega = \lambda \omega$ (lemma 1.4).
\qed \\
\subsection{Main theorem}
In this section we will prove the following theorem.
\begin{Thm}[Main theorem]\label{symres-main}
Let $\g$ be a semi-simple complex Lie algebra, and $G$ its adjoint group. Consider a nilpotent orbit $\0$ in $\g$, equipped with the Kostant-Kirillov 
form $\omega$. Then for any symplectic resolution
$\pi: (Z, \Omega) \rightarrow (\overline{\0}, \omega)$, there exists a parabolic subgroup $P$ of $G$, such that $(Z, \Omega)$ is isomorphic to 
$(T^*(G/P), \Omega_{can})$, where $\Omega_{can}$ is the canonical symplectic form on $T^*(G/P)$. Furthermore, under this isomorphism, the map $\pi$ becomes
$$T^*(G/P) \simeq G \times^P \mathfrak{u} \rightarrow \g, \qquad (g, X) \mapsto Ad(g)X ,$$ where $\mathfrak{u}$ is the nilradical of 
$\mathfrak{p} = Lie(P)$.  
\end{Thm}

The idea of the proof is to study the $\cit^*$-action on $Z$, as done in the papers of D. Kaledin \cite{Ka1} and H. Nakajima \cite{Nak}. 
A detailed account of the general theory can be found in section 6 of \cite{Ka1}. One may also consult the excellent book of H. Nakajima \cite{Na2}. 

By our proposition \ref{symI-lift}, the $\cit^*$-action on $\overline{\0}$ can be lifted on $Z$.  Let $Z^{\cit^*}$ be the fixed points 
subvariety in $Z$ under this
 $\cit^*$-action. Put $2n = dim(Z)$.

\begin{Lem}
There exists a $G$-equivariant attraction $p: Z \rightarrow Z^{\cit^*}$
\end{Lem}
\proof
For any point $x \in Z$, define $\phi_x: \cit^* \rightarrow Z$ to be $\phi_x(\lambda) = \lambda \cdot x$.  Let 
$\psi_x: \cit \rightarrow \overline{\0}$ be $\psi_x(\lambda) = \lambda \pi(x)$. Since $\pi$ is $\cit^*$-equivariant, we have the following
commutative diagram: \begin{equation*}
\begin{CD} \cit^* @>\phi_x>>  Z \\
@VVV      @VV{\pi}V \\
\cit @>\psi_x>> \overline{\0}
\end{CD} \end{equation*}
Applying the valuative criterion of properness for the proper map $\pi$, we  get a unique morphism $\widetilde{\phi}_x: \cit \rightarrow Z$ extending
 $\phi_x$. It is clear that $\widetilde{\phi}_x(0) \in Z^{\cit^*}$. 
 Now we define $p: Z \rightarrow Z^{\cit^*}$ by $p(x) = \widetilde{\phi}_x(0)$, which is the attraction map for the $\cit^*$-action. 
That the $G$-action commutes with the $\cit^*$-action on $Z$ implies that $p$ is $G$-equivariant.  
\qed \\

Now take a fixed point $z \in Z^{\cit^*}$, the action of $\cit^*$ on $Z$ induces a weight decomposition 
$$T_zZ = \oplus_{p \in \zit} T_z^p Z, $$ where $T_z^p Z = \{ v \in T_z Z |  \lambda_* v = \lambda^p v \}$. 
\begin{Def}
The $\cit^*$-action 
 is called {\em definite} at the point $z$ if $T_z^pZ = 0$ for all $p < 0$. 
\end{Def}
\begin{Lem}\label{symI-lem1}
There exists an irreducible smooth component $Z_0$ of $Z^{\cit^*}$ on which the $\cit^*$-action is definite.
\end{Lem}
\proof
Consider the open-dense set $\pi^{-1}(\0)$, which is isomorphic to $\0$, so $G$ acts on it transitively. By our preceding lemma, $p$ is $G$-equivariant,
so  the $G$-action on  $p(\pi^{-1}(\0))$ is also transitive. As a corollary, $p(\pi^{-1}(\0))$ is connected, then it is contained in a 
connected component, say $Z_0$, of $Z^{\cit^*}$. This means that the attraction subvariety of $Z_0$ contains the open-dense set $\pi^{-1}(\0)$ of $Z$, so
the $\cit^*$-action is definite on $Z_0$ (see lemma 6.1 \cite{Ka1}), i.e. $T_z^qZ=0$ for $q < 0$ and $\forall z \in Z_0$.
 Since $Z$ is smooth, it is well-known that the fixed points variety $Z^{\cit^*}$ is 
a union of smooth connected components, thus $Z_0$ is smooth and irreducible.
\qed \\
\begin{Lem}
The closed subvariety $Z_0$ is proper, $n$-dimensional and Lagrangian w.r.t. $\Omega$. 
\end{Lem}
\proof
Since $\pi$ is proper, the variety $\pi^{-1}(0)$ is proper. That $\pi$ is $\cit^*$-equivariant implies that
$Z_0$ is a closed subvariety in $\pi^{-1}(0)$, so $Z_0$ is proper.

The tangent space of $Z_0$ at the point $z$ equals to $T_z Z_0 = T_z^0 Z$. 
Now take two vectors $v_1 \in T_z^p Z$ and $v_2 \in T_z^qZ$, then the equation $\lambda^* \Omega = \lambda \Omega$
implies 
$$\lambda \Omega_z(v_1,v_2) = (\lambda^* \Omega)_z (v_1, v_2) = \Omega_{\lambda\cdot z} (\lambda_*(v_1), \lambda_*(v_2)) = \lambda^{p+q} \Omega_z(v_1,v_2).$$
So $\Omega_z(v_1,v_2)=0$ if $p+q \neq 1$. Now that $T_z^qZ=0$ for $q < 0$ implies that  for $p \geq 2$, the space $T_z^pZ$ is orthogonal to $T_zZ$ w.r.t. 
$\Omega$, thus $T_z^p Z = 0$, so $T_zZ = T_z^0Z \oplus T_z^1Z$. 
Furthermore $\Omega$ gives a duality between $T_z^0Z$ and $T_z^1Z$, so $dim(T_z^0Z) = n$, i.e. $Z_0$ is of dimension $n$. This also gives 
 that $Z_0$ is Lagrangian with 
respect to the symplectic form $\Omega$.
\qed \\

Since the decomposition by the attraction subvarieties is locally closed, the variety  $V = p^{-1}(Z_0)$ containing $\pi^{-1}(\0)$ is open-dense in $Z$.
As easily seen, it is also $G$-invariant. From now on, we will only consider $p: V \rightarrow Z_0$ instead of
 $p: Z \rightarrow Z^{\cit^*}$.  At the end of this section, we will prove that $V = Z$. 

 There exists a canonical symplectic 2-form $\Omega_{can}$ on  the cotangent bundle $T^*Z_0$, which comes from the Liouville form. There is also a natural
$\cit^*$-action on $T^*Z_0$, considered as a vector bundle over $Z_0$. 
\begin{Lem}
There exists a $\cit^*$-equivariant isomorphism $i: V \rightarrow T^*Z_0$, which identifies also the two symplectic structures.
\end{Lem}
\proof
Since the $\cit^*$-action is definite on $Z_0$,
a classical work of Bialynicki-Birula (see Theorem 2.5 \cite{BB}) implies  that the attraction  $p: V \rightarrow Z_0$
is a locally trivial (for the Zariski topology) fibration with vector spaces as fibers. Furthermore, the decomposition $T_zZ = T_z^0Z \oplus T_z^1Z$
shows that $\cit^*$ acts on fibers linearly, thus  the map $p$ makes $V$ a vector bundle of rank $n$ 
over $Z_0$, and the $\cit^*$-action on this vector bundle is the natural one. 
Let us identify $Z_0$ with the zero section of this bundle.  Since for any $z \in Z_0, T_zZ = T_zZ_0 \oplus T^1_zZ$, the induced $\cit^*$-action on
 the normal bundle $N$ of $Z_0$ in $V$ is the natural one when we consider
$N$ as a vector bundle over $Z_0$. This gives a $\cit^*$-equivariant  isomorphism between $V$ and the total space of the normal bundle $N$. Let us also 
denote by $\Omega$ the symplectic form on $N$, which comes from the symplectic form on $V$.

Now we establish an isomorphism between $(N,\Omega)$ and $(T^*Z_0, \Omega_{can})$ as follows. Take a point $z \in Z_0$, and a vector
$v \in N_z$. Since $Z_0$ is Lagrangian in the two symplectic spaces, there exists a unique vector $w \in T_z^*Z_0$ such that
 $\Omega_z(v,u) = \Omega_{can,z}(w,u)$ for all $u \in T_zZ_0$.
We define the map $i: N \rightarrow T^*Z_0$ to be $i(v)=w$. Now it is easy to see that $i$ is a symplectic isomorphism.
\qed \\

From now on, we will denote by $\Omega$ (instead of $\Omega_{can}$) the canonical symplectic 2-form on $T^*Z_0$.
Now we will study the action of $G$ on $Z$. 
Since the $G$-action commutes with the $\cit^*$-action, the fixed points $Z^{\cit^*}$ is $G$-invariant. Since $Z_0$ is a connected component of 
$Z^{\cit^*}$ and $G$ is connected,  $Z_0$ is also $G$-invariant, thus we get a $G$-action on $Z_0$. 
This action induces naturally an action of $G$ on the total space of $T Z_0$, which is given by $g\cdot(z,v) = (g\cdot z, dg(v))$.
 By taking the dual, we have a
$G$-action on $T^*Z_0$. 
\begin{Lem}\label{symI-lem2}
The $G$-action on $T^*Z_0$ is isomorphic to the one described above. 
\end{Lem}
\proof
For an element $g \in G$, we write $\phi_g$ the action of $g$ on $T^*Z_0$. The action of $G$ on $Z_0$ will be  simply denoted by $\cdot$.
For an element $z \in Z_0 \subset T^*Z_0$, we have a natural decomposition $T_z(T^*Z_0) = T^*_zZ_0 \oplus T_zZ_0$, which is also isotropic w.r.t. $\Omega$.
The $G$-action on $T_zZ_0$ is the one induced by the action of $G$ on $Z_0$.
Since the $G$-action commutes with the $\cit^*$-action, for a vector $v_2 \in T_z^*Z_0$, considered as a tangent  vector in $T_{z}(T^*Z_0)$,   
we have $d\phi_g(v_2) = \frac{d}{d\lambda}|_{\lambda = 0} \phi_g(\lambda v_2) = \frac{d}{d\lambda}|_{\lambda = 0} \lambda \phi_g( v_2)= \phi_g(v_2)$. 

For any vector  $v_1 \in T_zZ_0$, since $\Omega$ is $G$-invariant,  we have   
$$ \Omega_z(v_1,v_2) = \phi_g^*\Omega_z(v_1, v_2) = \Omega_{g\cdot z}(d\phi_g(v_1), d\phi_g(v_2)) = \Omega_{g\cdot z}(dg(v_1),\phi_g(v_2)).$$
That $\Omega$ is the natural symplectic structure on $T^*Z_0$ implies that  $\phi_g(v_2) = (dg)^*(v_2)$, i.e. the action of $G$ on $T^*Z_0$ is the natural one 
induced from the action of $G$ on $Z_0$.
\qed \\
\begin{Lem}\label{symI-lem3}
The variety $Z_0$ contains an open-dense $G$-orbit.
\end{Lem}
\proof
The argument in the proof of lemma \ref{symI-lem1} shows that $G$ acts transitively on the set $p(\pi^{-1}(\0)) \subset Z_0$.
Here $p: T^*Z_0 \rightarrow Z_0$ is the canonical projection. Since $\pi^{-1}(\0)$ is open, $ p(\pi^{-1}(\0))$ is an open-dense $G$-orbit. 
\qed \\

The author wants to thank M. Brion and  F. Knop  for having provided the proof of the following proposition, which is implicit in \cite{Kno}.
\begin{Prop}\label{symI-km}
Let $\g$ be a semi-simple complex Lie algebra and $G$ its adjoint group. Let $P$ be a closed subgroup of $G$ and $\mathfrak{p}$ its Lie algebra.
Then $P$ has a dense orbit (via the co-adjoint action) in $(\g/\mathfrak{p})^*$ if and only if $P$ is parabolic. 
\end{Prop}
\proof
Suppose that $P$ has a dense orbit in $(\g/\mathfrak{p})^*$, then its image in $\g^*//G$ is a point. This means that the space $G/P$ has rank 0
(see Satz 5.4 \cite{Kno}).  Now Satz 9.1 of {\em loc. cit.} implies that all isotropy groups of $G/P$ are parabolic, thus $P$ is parabolic.

The converse is a theorem of Richardson \cite{Ric}.
\qed \\
\begin{Cor}
The $G$-action on $Z_0$ is transitive, and $Z_0$ is isomorphic to $G/P$ for some parabolic subgroup $P$ of $G$.
\end{Cor}
\proof
By our lemma \ref{symI-lem3}, $Z_0$ has an open $G$-orbit, say $U = G/P$, where $P$ is a subgroup of $G$. Now $T^*(G/P) \simeq G \times^P (\g/\mathfrak{p})^*$.
Note that $G$ has an open orbit in $T^*(G/P)$, which is isomorphic to $\pi^{-1}(\0)$, so $P$ has an open dense orbit in $(\g/\mathfrak{p})^*$, where
the action is the co-adjoint one by our lemma \ref{symI-lem2}. Now the above proposition gives that $P$ is
parabolic. As a corollary, the variety $G/P$ is projective which is also dense in the proper irreducible variety $Z_0$, so $Z_0 = G/P$.
\qed \\

We are indebted to M. Brion for providing ideas to the proof of the following 
\begin{Prop} \label{orbitcover}
Let $\0, \0'$ be two nilpotent orbits in a semisimple complex Lie algebra $\g$.
Let $G$ be the adjoint group of $\g$. 
Suppose that there exists a finite $G$-equivariant morphism $\pi: \0 \to \0'$, which commutes with 
the $\cit^*$-action.
Then $\0 = \0'$.
\end{Prop}
\proof
 Let $X = P(\overline{\0})$ and $X' = P(\overline{\0'})$. Then we have finite  surjective $G$-morphisms 
$P(\widetilde{\0}) \to P(\widetilde{\0'})$, $P(\widetilde{\0}) \to X$ and $P(\widetilde{\0'}) \to X'$,
This implies that the moment polytopes (for definition see \cite{Bri})
 associated to these varieties are the same, which will be denoted by $\C$.

 Let $\beta$ be the point closest to 0 in $\C$. Consider  the set $S_\beta$  of
points $\bar{x} \in X$ such that $\beta$ is the closest point to 0 in $\C(\overline{G \cdot \bar{x}})$.
Then $P(\0) \subset S_\beta$. By Proposition 3.2 [Bri], there exists a $G$-morphism
$\phi: S_\beta \to C_{n \beta}$. Notice that the $G$-action on $C_{n \beta}$ is
transitive, so there exists a point $\bar{x} \in  P(\0)$ such that $\phi(\bar{x}) = x_{n \beta}$.
Now by Proposition 3.3 [Bri], the one-parameter subgroup $\beta / q(\beta)$ is optimal for 
$x$, where $x \in \0$ is a pre-image of $\bar{x}$.

Similarly we find a point $x' \in \0'$ which has $\beta / q(\beta)$ as an optimal 
one-parameter subgroup. Now by Cor. 11.6 \cite{Hes2}, $x$ is conjugate to $x'$, thus $\0 = \0'$.
\qed \\

By our above lemmas,  $V$ is isomorphic to $T^*(G/P) = G \times^P \mathfrak{u}$, where $\mathfrak{u}$ is the nilradical of $\mathfrak{p} = Lie(P)$ which is
 identified with $(\g/\mathfrak{p})^*$ via the Killing form.
Under this isomorphism, the $G$-action on $G \times^P \mathfrak{u}$ is the natural one $g \cdot (h, X) = (g h, X)$.
\begin{Lem} \label{gap}
Under these isomorphisms, the map $\pi: V \simeq G \times^P \mathfrak{u}  \rightarrow \overline{\0}$ becomes $(g,X) \mapsto Ad(g)X$. 
\end{Lem}
\proof
Take a Richardson element $X \in \mathfrak{u}$, then $G \cdot (1,X)$ is dense in $V = G \times^P \mathfrak{u}$, 
thus if we let $Y = \pi(1,X)$, then $G \cdot Y$ is dense in $\overline{\0}$, which shows that $Y \in \0$. 
From that $\pi$ is a resolution, we have an isomorphism between $\0 = G \cdot Y$ and 
$G \cdot (1,X) \simeq G/P^X$. Let $\0_X$ be the nilpotent orbit of $X$, then we obtain an unramified finite $G$-covering from $\0$ to $\0_X$, 
which commutes with the $\cit^*$-action. By Proposition \ref{orbitcover},  $\0 = \0_X$, i.e. the two nilpotent orbits coincide.
Thus $X \in \0$, and there exists an element $g \in G$ such that $g Y = X$. By using the action of $g$ on $G/P$, 
 we can suppose $\pi(1,X) = X$.
Since $\pi$ is $G$-equivariant, for any $p \in P$, we have $\pi(1,Ad(p)X) = \pi(1, p\cdot X) = p \cdot \pi(1,X) = Ad(p)X$.
 Now that $P \cdot X$ is dense in $\mathfrak{u}$ implies that $\pi(1,Y) = Y$ for any $Y \in \mathfrak{u}$, which gives   
 $\pi(g, Y) = \pi(g \cdot(1,Y)) = g \cdot \pi(1,Y) = Ad(g)Y$ for any $g \in G$.
\qed \\

Now by [BK], the map $\pi: V \to \overline{\0}$ is proper. Thus the properness of the map
$Z \to \overline{\0}$ implies that $V= Z$,  which concludes the proof of the Main theorem. 

\subsection{Symplectic resolutions: classical cases}
Let $\g$ be a simple Lie algebra and $G$ its adjoint group. Recall that a nilpotent orbit $\0$ in $\g$ is called 
{\em polarizable} if there exists a parabolic  subalgebra $\p \subseteq \g$ such that $\mathfrak{u} \cap \0$ is dense in $\mathfrak{u}$, where 
$\mathfrak{u}$ is the nilradical of $\p$. The Lie group $P$ corresponding to $\p$ is called a {\em polarization} of $\0$. In the literature, these orbits are 
also called  {\em Richardson orbits}. For such an orbit, we have $dim(\0) = 2\, dim(G/P)$. As shown by Richardson, every parabolic subalgebra corresponds 
to a polarizable orbit. 
In the case of $\g = \mathfrak{sl}_n$, every nilpotent orbit is polarizable. As a direct corollary of our main theorem, we have
\begin{Prop}
If the symplectic variety $\overline{\0}$ admits a symplectic resolution, then $\0$ is a Richardson orbit.
\end{Prop}

So in the following we will only consider Richardson orbits.
For an element $X \in \0$, let $G^X$ be the stabilizer of $X$ in $G$, and $(G^X)^\circ$ the identity component. The component group $G^X/(G^X)^\circ$
is denoted by $A(\0)$. Suppose that $\0$ is polarizable and let $P$ be a polarization. Let $P^X$ be the stabilizer of $X$ in $P$. Then we have 
$(G^X)^\circ \subseteq P^X \subseteq G^X$. Put $A_P(\0) = P^X/(G^X)^\circ$, which is also the stabilizer of $P$ in $A(\0)$. Note that $A_P(\0)$ is a 
subgroup  of the  group $A(\0)$. By Corollary 6.1.7 \cite{CM}, $A(\0)$ is abelian.
 We denote by $N(P) = [A(\0) : A_P(\0)]$ the index of $A_P(\0)$ in $A(\0)$.
It turns out that this number is useful to our problem, as shown by the following:
\begin{Prop}
Let $\0$ be a nilpotent orbit in a simple Lie algebra $\g$.  Then the symplectic
variety $\overline{\0}$ admits a symplectic resolution if and only if 
there exists a polarization $P$ of $\0$ such that $N(P) = 1$.
\end{Prop}
\proof
By our main theorem, any symplectic resolution for $\overline{\0}$ is of the form
 $$\pi: T^*(G/P) \simeq G \times^P \mathfrak{u} \rightarrow \overline{\0}, \quad \pi(g,X) = Ad(g)X,  $$
for some parobolic subgroup $P$ in $G$. In particular, we see that $P$ gives a polarization of $\0$. 
As shown in \cite{BK}(\S 7), the map $\pi$ is proper and of degree $N(P) = [G^X:P^X]$. 
So it is birational if and only if $N(P) = 1$, in which case it 
  gives an isomorphism  from $\pi^{-1}(\0)$ to $\0$, thus we get a symplectic resolution for $\overline{\0}$.
\qed \\
\begin{Cor}\label{symreso-slnil}
The closure of any nilpotent orbit in $\mathfrak{sl}_n$ admits a symplectic resolution.
\end{Cor}
This corollary comes directly from the proposition and the fact that $A(\0) =1$ for any nilpotent orbit in $\mathfrak{sl}_n$ (see for example corollary
6.1.6 of \cite{CM}). The following corollary follows from corollary \ref{symI-cor}.
\begin{Cor}
Consider a nilpotent orbit $\0$ in $\mathfrak{sl}_n$. If $\overline{\0} - \0$ is not of pure codimension 2, then  $\overline{\0}$ is not locally
$\qit$-factorial.
\end{Cor}

In the following we will consider the classical cases $\g = \mathfrak{so}_m$ and $\g =\mathfrak{sp}_m$. The first question is to decide when a nilpotent 
orbit $\0$ is polarizable, and the second question is to decide when there exists a polarization $P$ such that $N(P) =1$. We will use notations and results 
of \cite{Hes} to settle these two questions. 

First, some notations. Here all congruences are modulo 2. Let $\varepsilon = 0$ for $\mathfrak{so}_m$ and $\varepsilon = 1$ for $\mathfrak{sp}_m$. A 
natural number $q \geq 0$ is called {\em admissible} if $q \equiv m$ and $q \neq 2$ if $\varepsilon = 0$. Set 
$$Pai(m,q) = \{partitions \,\pi \,of\, m | \, \pi_j \equiv 1 \, if \, j \leq q;\, \pi_j \equiv 0 \, if \, j > q \}, $$
 which parametrizes some Levi types of parabolic subgroups in $G$.  Let $\mathcal{P}_\varepsilon(m)$ be the partitional parameterizations of 
nilpotent orbits in $\g$, i.e. the partitions {\bf d} such that for any $i \equiv \varepsilon$, the number $\#\{j|d_j = i\}$ is even 
(see section 5.1 in \cite{CM}). 
 Then there exists an injective Spaltensein mapping $$ S_q: Pai(m,q) \rightarrow \mathcal{P}_\varepsilon(m). $$

As shown by theorem 7.1.(a) in \cite{Hes}, a nilpotent orbit $\0$ is polarizable if and only if there exists some admissible $q$ such that the partition of
 the orbit is in the image of $S_q$. If we denote by $$J(d) = \{j | d_j \equiv \varepsilon \} \cup \{j, j+1 | j \equiv m, d_j = d_{j+1} \} $$
$$j_1(d) = sup \{ j \in J(d) | d_j \equiv 1 \} \ ( - \infty \ if\ empty )$$ $$ j_0(d) = min\{ j \in J(d) | d_j \equiv 0 \} \ ( + \infty\ if \ empty),$$
then a partition {\bf d} is in the image of $S_q$ if and only if (see Prop. 6.5 \cite{Hes}): 
$$j_1(d) \leq q <j_0(d) \ and \  d_j \equiv d_{j+1} \ if\ j \equiv m+1. $$
Let $B(d) = \{ j \in \nit | d_j > d_{j+1}; d_j \equiv \varepsilon+1 \}$ and $u = \frac{1}{2}(-1)^\varepsilon ( \#\{j | d_j \equiv 1\} - q )$, then we have 
\begin{Prop}[Theorem 7.1 \cite{Hes}]
Suppose {\bf d}  is in the image of $S_q$. Let $P$ be an associated  parabolic subgroup. Then 
$$N(P) = \begin{cases} & 2^u \quad if \ q+\varepsilon \geq 1 \ or \ B(d) = \emptyset; \\  & 2^{u-1} \quad if \ q=\varepsilon=0 \ and \ B(d)
 \neq \emptyset.  \end{cases} $$
\end{Prop}
Now we will do a case-by-case check. 
\begin{Prop}\label{symreso-sp}
For $\g = \mathfrak{sp}_{2n}$ (resp. $\g = \mathfrak{so}_{2n+1}$), let $\0$ be a nilpotent orbit in $\g$ corresponding to the partition 
{\bf d}= $[d_1, \cdots, d_N]$. Then the following three conditions are equivalent:

(i) $\0$ is polarizable and there exists a polarization $P$ such that $N(P) = 1$;

(ii) there exists an even (resp. odd) number $q \geq 0$ such that the first $q$ parts $d_1, \cdots, d_q$ are odd and the other parts are even;

(iii) the symplectic variety $\overline{\0}$  has a symplectic resolution.
\end{Prop}
\proof
The equivalence between (i) and (iii) has been established earlier. In the following, we will prove the equivalence between conditions (i) and  (ii).
 
For the case $\g =\mathfrak{sp}_{2n}$, we have $\varepsilon = 1$, so $N(P) = 2^u$, which gives that $N(P)=1$ if and only if $u=0$, i.e. if and only if
$q = \# \{j | d_j \equiv 1 \}$. Now the condition (i) is equivalent to saying that the partition {\bf d} is in the image of 
$S_q: Pai(2n,q) \rightarrow \mathcal{P}_1(2n)$. This condition is equivalent to having  $j_1(d) \leq q < j_0(d)$ and $d_j \equiv d_{j+1}$ if $j \equiv 1$. 
Note that here the set $J(d)$ contains the set $\{j | d_j \equiv 1\}$, thus $j_1(d) = sup \{j | d_j \equiv 1 \}$. Now
 $j_1(d) \leq q = \# \{j | d_j \equiv 1 \}$ implies that the first $q$ parts of {\bf d} are odd, and the others are even.  In this case, the conditions
$q < j_0(d)$ and  $d_j \equiv d_{j+1}$ if $j \equiv 1$ are satisfied automatically. So we have the equivalence of the conditions (i) and (ii).

In the case of $\g = \mathfrak{so}_{2n+1}$, since $q$ is admissible,  $q \equiv 1$, thus $q \geq 1$, which gives $N(P) = 2^u$. So $N(P) = 1 $ if and 
only if $u=0$, i.e. $q = \# \{j | d_j \equiv 1 \}$. Now the set $\{j | d_j \equiv 0\}$ is contained in $J(d)$, thus $j_0(d) = min\{j | d_j \equiv 0\}$.
The condition $j_0(d) > q$ gives that the first $q$ parts should be odd, and the others are even. If this is satisfied, then the orbit is polarizable
by some $P$ with $N(P) = 1$. So we have the equivalence between (i) and (ii).
\qed \\

{\em Examples:} Consider $\mathfrak{so}_7$, the variety $\overline{\0}_{[3,2,2]}$ has a symplectic resolution. So in the case of $\mathfrak{so}_7$, only 
$\overline{\0}_{min}$ does not admit any symplectic resolution.
For $\mathfrak{sp}_6$, we see that the variety  $\widetilde{\0}_{[4,1,1]}$ does not admit any symplectic resolution. 

\begin{Prop}\label{symreso-so}
For $\g = \mathfrak{so}_{2n}$, let $\0$ be a nilpotent orbit corresponding to the partition {\bf d}=$[d_1, \cdots, d_N]$. Then the following three
conditions are equivalent:

(i) $\0$ is polarizable and there exists a polarization $P$ such that $N(P) = 1$;

(ii) either there exists some even number $q \neq 2$ such that the first $q$ parts of {\bf d} are odd and the others are even or there exists 
exactly 2 odd parts which are at the positions $2k-1$ and $2k$ in the partition {\bf d} for some $k$;

(iii) the symplectic variety $\overline{\0}$ has a symplectic resolution.
\end{Prop}
\proof
We need to establish the equivalence between (i) and (ii). There are three cases:

{\em Case (a):} $q =\varepsilon = 0$ and $B(d) \neq \emptyset$. 

In this case we should have $u = 1$, then $\# \{j | d_j \equiv 1 \} = 2$,
 i.e. there are exactly 2 odd parts. The polarizable condition $d_j \equiv d_{j+1}$ if $j \equiv 1$ gives that the two parts are at the positions
$2k-1$ and $2k$ in {\bf d} for some $k$. In this case, the other conditions as $B(d) \neq \emptyset$ and $j_1(d) \leq 0 < j_0(d)$ are satisfied.

{\em Case (b):} $B(d) = \emptyset$

In this case we find that there is no odd part in {\bf d}.

{\em Case (c):} $q \geq 2$ and $q \equiv 0$.

Since $q$ should be admissible, thus $q \neq 2$, i.e. $q \geq 4$ even. Now $N(P) = 1$ gives $u = 0$, i.e. $q = \# \{j | d_j \equiv 1 \}$. 
Now the set $\{j | d_j \equiv 0\} $ is contained in $J(d)$, so $j_0(d) > q$ gives that the first $q$ parts of the partition {\bf d}
should be odd, and thus the other parts should be even. If this is satisfied, then the other conditions are also satisfied.

The above analysis for the three cases gives the equivalence between the conditions (i) and (ii). 
\qed \\

\subsection{Symplectic resolutions: exceptional cases}
Firstly we have the following
\begin{Prop}\label{symI-excep}
Let $\g$ be one of the following exceptional simple complex Lie algebras: 
$G_2, F_4, E_6$, and $\0$ a nilpotent orbit in $\g$. Then the symplectic variety
$\overline{\0}$ admits a symplectic resolution if and only if $\0$ is a Richardson orbit.
\end{Prop}
\proof
Note that an even orbit is always a Richardson orbit, and it admits a symplectic resolution by Springer's resolution.
Among the non-even orbits,  a list of Richardson orbits in exceptional simple Lie algebra $\g$ can be found in \cite{Hir}.
Here we use notations of \cite{CM} for nilpotent orbits in exceptional algebras.
When $\g = G_2$, every Richardson orbit is even. When $\g = F_4$, there is only one non-even Richardson orbit: $C_3$. In the case of
$\g = E_6$, there are 5 non-even Richardson orbits: $2A_1, A_2+2A_1, A_3, A_4+A_1, D_5(a_1)$. Now the table given in \cite{CM} (chap. 8)
shows that  all these orbits are simply connected, thus $A(\0) = 1$, so for any polarization $P$ for $\0$, the collapsing $T^*(G/P) \rightarrow \overline{\0}$
gives a symplectic resolution. 
\qed \\

In the case of $\g = E_7$, there are 5 non-even Richardson orbits, two of which have trivial component groups (thus admit a symplectic resolution):
$D_5+A_1$ and $D_6(a_1)$; while for the other three ($D_4(a_1)+A_1, A_4+A_1$ and $D_5(a_1)$) have component groups  $S_2$, so the above argument does
not apply. We don't know whether they admit a symplectic resolution. 

In the case of $\g = E_8$, there are 7 non-even Richardson orbits, three of which have trivial component groups (thus admit a symplectic resolution):
$A_4+A_2+A_1, A_6+A_1$ and $E_7(a_1)$; while for the other four ($D_6(a_1), D_7(a_2), E_6(a_1)  + A_1$ and $E_7(a_3)$), whose component groups are $S_2$.
As above, we don't know the answer.

\quad \\
Labortoire J.A.Dieudonn\'e, Parc Valrose \\ 06108 Nice cedex 02, FRANCE \\
baohua.fu@polytechnique.org
\end{document}